\documentclass[a4paper,12pt]{article}
\usepackage{enumerate}

\makeatletter

\@setfontsize\normalsize\@xiipt{15.5}%

\renewcommand*{\author}[1]{\gdef\@author{#1}\gdef\@pauthor{{\def\and{ --- }#1}}}
\renewcommand*{\title}[1]{\gdef\@title{#1}\gdef\@ptitle{#1}}
\if@twoside         % If two-sided printing.
  \def\ps@draft{%
    \def\@oddfoot{\small\null\hfil\thepage\hfil}% Page number in feet
    \let\@evenfoot\@oddfoot
    \def\@evenhead{\small\@date\hfil\slshape\@pauthor\hfil}% Left head
    \def\@oddhead{\small\null\hfil\slshape\@ptitle\hfil}% Right head
    \let\@mkboth\@gobbletwo% Disable marking commands
    \let\sectionmark\@gobble
    \let\subsectionmark\@gobble
   }
\else               % If one-sided printing.
  \def\ps@draft{%
    \def\@oddfoot{\small\@date\hfil\slshape\@pauthor\hfil\upshape\thepage}
    \def\@oddhead{\small\null\hfil\slshape\@ptitle\hfil}
    \let\@mkboth\@gobbletwo% Disable marking commands
    \let\sectionmark\@gobble
    \let\subsectionmark\@gobble
  }
\fi

\newcommand*{\keywords}[1]{\gdef\@keywords{#1}}
\newcommand{\keywordsname}{Key words and phrases}

\newcommand*{\subjclass}[1]{\gdef\@subjclass{#1}}
\newcommand{\subjclassname}{1991 AMS Mathematics Subject Classification}

\def\@maketitle{%
  \newpage
  \null
  \vskip 2em%
  \begin{center}%
    {\Large\bfseries \@title \par}%
    \vskip 1.5em%
    {\small\scshape
      \lineskip .5em%
      \begin{tabular}[t]{c}%
        \@author
      \end{tabular}\par
    }%
    \vskip 1em%
    {\small\@date}
  \end{center}%
  \par
  \vskip 1.5em
  \begingroup
    \let\@makefnmark\relax \let\@thefnmark\relax
    \ifx\@empty\@subjclass\else
       \@footnotetext{{\itshape\subjclassname}.\enspace\@subjclass.}
    \fi
    \ifx\@empty\@keywords\else
       \@footnotetext{{\itshape\keywordsname}.\enspace\@keywords.}
    \fi
  \endgroup
}

\renewcommand{\part}{\par
   \addvspace{4ex}%
   \@afterindentfalse
   \secdef\@part\@spart}

\def\@part[#1]#2{%
    \ifnum \c@secnumdepth >\m@ne
      \refstepcounter{part}%
      \addcontentsline{toc}{part}{\thepart\hspace{1em}#1}%
    \else
      \addcontentsline{toc}{part}{#1}%
    \fi
    {\parindent \z@ \raggedright
     \interlinepenalty \@M
     \reset@font
     \ifnum \c@secnumdepth >\m@ne
       \large\bfseries \partname~\thepart
       \par\nobreak
     \fi
     \Large \bfseries #2%
     \markboth{}{}\par}%
    \nobreak
    \vskip 3ex
    \@afterheading}
\def\@spart#1{%
    {\parindent \z@ \raggedright
     \interlinepenalty \@M
     \reset@font
     \Large \bfseries #1\par}%
     \nobreak
     \vskip 3ex
     \@afterheading}
\def\@endpart{\vfil\newpage
              \if@twoside
                \hbox{}%
                \thispagestyle{empty}%
                \newpage
              \fi
              \if@tempswa
                \twocolumn
              \fi}
\renewcommand{\section}{\@startsection {section}{1}{\z@}%
                                   {-3.5ex \@plus -1ex \@minus -.2ex}%
                                   {2.3ex \@plus.2ex}%
                                   {\reset@font\large\bfseries}}
\renewcommand{\subsection}{\@startsection{subsection}{2}{\z@}%
                                     {-3.25ex\@plus -1ex \@minus -.2ex}%
                                     {1.5ex \@plus .2ex}%
                                     {\reset@font\normalsize\bfseries}}
\renewcommand{\subsubsection}{\@startsection{subsubsection}{3}{\z@}%
                                     {-3.25ex\@plus -1ex \@minus -.2ex}%
                                     {1.5ex \@plus .2ex}%
                                     {\reset@font\normalsize\bfseries}}
\renewcommand{\paragraph}{\@startsection{paragraph}{4}{\z@}%
                                    {3.25ex \@plus1ex \@minus.2ex}%
                                    {-1em}%
                                    {\reset@font\normalsize\bfseries}}
\renewcommand{\subparagraph}{\@startsection{subparagraph}{5}{\parindent}%
                                       {3.25ex \@plus1ex \@minus .2ex}%
                                       {-1em}%
                                      {\reset@font\normalsize\bfseries}}

\renewcommand{\theenumi}{\alph{enumi}}

\renewcommand{\theenumii}{\roman{enumii}}

\renewcommand{\p@enumii}{\theenumi.}
\renewcommand{\theenumiii}{\Alph{enumiii}}

\renewcommand{\p@enumiii}{\theenumi.\theenumii.}
\renewcommand{\p@enumiv}{\p@enumiii\theenumiii.}

\@addtoreset{equation}{section}

\RequirePackage{amsthm}

\renewenvironment{proof}[1][\proofname]{\par
  \normalfont
  \topsep6\p@\@plus6\p@ \trivlist
  \item[\hskip\labelsep\slshape
    #1\@addpunct{.}]\ignorespaces
}{%
  \qed\endtrivlist
}

\theoremstyle{plain}

\newtheorem{theorem}{Theorem}[section]
\newtheorem{proposition}[theorem]{Proposition}
\newtheorem{lemma}[theorem]{Lemma}
\newtheorem{corollary}[theorem]{Corollary}

\theoremstyle{definition}

\newtheorem{definition}[theorem]{Definition}
\newtheorem{notation}[theorem]{Notation}
\newtheorem{example}[theorem]{Example}

\newtheorem{remark}[theorem]{Remark}

\RequirePackage{amsmath}
\RequirePackage{amssymb}

\newcommand{\C}{\mathbb{C}}

\newcommand{\R}{\mathbb{R}}

\newcommand{\Z}{\mathbb{Z}}

\DeclareMathOperator{\id}{id}

\newcommand{\DSum}{\bigoplus}

\renewcommand{\to}[1][]{\xrightarrow[#1]{}}

\newcommand{\isoto}[1][]{\xrightarrow[#1]{\sim}}

\newcommand{\Endo}[1][]{\mathrm{End}_{\raise1.5ex\hbox to.1em{}#1}}

\newcommand{\Hom}[1][]{\mathrm{Hom}_{\raise1.5ex\hbox to.1em{}#1}}

\newcommand{\RHom}[1][]{\mathrm{RHom}_{\raise1.5ex\hbox to.1em{}#1}}

\newcommand{\Ext}[2][]{\mathrm{Ext}_{\raise1.5ex\hbox to.1em{}#1}^{#2}}

\newcommand{\THom}[1][]{\mathrm{THom}_{\raise1.5ex\hbox to.1em{}#1}}

\newcommand{\Tens}[1][]{\mathbin{\otimes_{\raise1.5ex\hbox to-.1em{}#1}}}

\newcommand{\LTens}[1][]{\mathbin{\otimes_{\raise1.5ex\hbox to-.1em{}#1}^{L}}}

\newcommand{\Tor}[2][]{\mathrm{Tor}^{\raise1.5ex\hbox to.1em{}#1}_{#2}}

\def\sha{\mathcal{A}}
\def\shb{\mathcal{B}}
\def\shc{\mathcal{C}}
\def\shd{\mathcal{D}}

\def\shf{\mathcal{F}}

\def\shi{\mathcal{I}}

\def\shl{\mathcal{L}}
\def\shm{\mathcal{M}}
\def\shn{\mathcal{N}}
\def\sho{\mathcal{O}}

\def\shu{\mathcal{U}}
\def\shv{\mathcal{V}}

\newcommand{\sect}{\Gamma}

\renewcommand{\hom}[1][]{{\mathcal{H}om}_{\raise1.5ex\hbox to.1em{}#1}}

\newcommand{\rhom}[1][]{{R\mathcal{H}om}_{\raise1.5ex\hbox to.1em{}#1}}

\newcommand{\ext}[2][]{{\mathcal{E}xt}_{\raise1.5ex\hbox to.1em{}#1}^{#2}}

\newcommand{\thom}[1][]{{T\mathcal{H}om}_{\raise1.5ex\hbox to.1em{}#1}}

\newcommand{\tens}[1][]{\mathbin{\otimes_{\raise1.5ex\hbox to-.1em{}#1}}}

\newcommand{\ltens}[1][]{\mathbin{\otimes_{\raise1.5ex\hbox to-.1em{}#1}^{L}}}

\newcommand{\tor}[2][]{{\mathcal{T}or}^{\raise1.5ex\hbox to.1em{}#1}_{#2}}

\DeclareMathOperator{\supp}{supp}

\newcommand{\oim}[1]{{#1}_*}

\newcommand{\opb}[1]{#1^{-1}}

\newcommand{\GHom}[1][]{\mathrm{GHom}_{\raise1.5ex\hbox to.1em{}#1}}

\newcommand{\GExt}[2][]{\mathrm{GExt}_{\raise1.5ex\hbox to.1em{}#1}^{#2}}

\newcommand{\FHom}[1][]{\mathrm{FHom}_{\raise1.5ex\hbox to.1em{}#1}}

\newcommand{\ghom}[1][]{{\mathcal{GH}om}_{\raise1.5ex\hbox to.1em{}#1}}

\newcommand{\gext}[2][]{{\mathcal{GE}xt}_{\raise1.5ex\hbox to.1em{}#1}^{#2}}

\newcommand{\fhom}[1][]{{\mathcal{FH}om}_{\raise1.5ex\hbox to.1em{}#1}}

\newcommand{\gr}{\mathop{\mathcal{G}r}\nolimits}

\newcommand{\tenstop}[1][]{\mathbin{\hat{\otimes}_{\raise1.5ex\hbox to-.1em{}#1}}}

\newcommand{\homtop}[1][]{\mathcal{L}_{\raise1.5ex\hbox to.1em{}#1}}

\newcommand{\Homtop}[1][]{\mathrm{L}_{\raise1.5ex\hbox to.1em{}#1}}

\newcommand{\D}{\mathcal{D}}

\newcommand{\E}{\mathcal{E}}

\newcommand{\ER}{\mathcal{E}^{\R}}

\def\absdoim#1{\underline{#1}_*}
\def\reldoim[#1]#2{\underline{#2}_{|{#1}*}}
\def\doim{\@ifnextchar [{\reldoim}{\absdoim}}

\def\absdeim#1{\underline{#1}_*}
\def\reldeim[#1]#2{\underline{#2}_{|{#1}*}}
\def\deim{\@ifnextchar [{\reldeim}{\absdeim}}

\def\absdopb#1{\underline{#1}^{-1}}
\def\reldopb[#1]#2{\underline{#2}_{|{#1}}^{-1}}
\def\dopb{\@ifnextchar [{\reldopb}{\absdopb}}

\def\absboim#1{\underline{\underline{#1}}_*}
\def\relboim[#1]#2{\underline{\underline{#2}}_{|{#1}*}}
\def\boim{\@ifnextchar [{\relboim}{\absboim}}

\def\absbeim#1{\underline{\underline{#1}}_*}
\def\relbeim[#1]#2{\underline{\underline{#2}}_{|{#1}*}}
\def\beim{\@ifnextchar [{\relbeim}{\absbeim}}

\def\absbopb#1{\underline{\underline{#1}}^*}
\def\relbopb[#1]#2{\underline{\underline{#2}}_{|{#1}}^*}
\def\bopb{\@ifnextchar [{\relbopb}{\absbopb}}

\RequirePackage{xypic}

\newcommand{\contatto}{\alpha}
\newcommand{\simplettico}{\omega}

\newcommand{\OO}{\sho}
\newcommand{\field}{\mathbf{k}}

\newcommand{\catMod}{\mathrm{Mod}}
\newcommand{\atw}[1]{{}_{#1}}

\renewcommand{\D}[1][]{\mathcal{D}_{#1}}

\newcommand{\halfform}{{\sqrt v}}

\renewcommand{\E}[1][]{\mathcal{E}_{#1}}
\renewcommand{\ER}[1][]{\mathcal{E}_{#1}^{\R}}
\newcommand{\Ev}[1][]{\mathcal{E}_{#1}^\halfform}
\newcommand{\ERv}[1][]{\mathcal{E}_{#1}^{\halfform,\R}}
\newcommand{\Iv}[1][]{\mathcal{I}_{#1}^\halfform}
\newcommand{\Lrm}{\mathrm{L}}

\newcommand{\EE}[1][]{\mathsf{E}_{#1}}
\newcommand{\EER}[1][]{\mathsf{E}_{#1}^{\R}}

\newcommand{\W}[1][]{\mathcal{W}_{#1}}
\newcommand{\Wv}[1][]{\mathcal{W}_{#1}^\halfform}

\newcommand{\WW}[1][]{\mathsf{W}_{#1}}

\renewcommand{\gr}{\mathop{\mathrm{gr}}}
\renewcommand{\DSum}{\bigoplus}
\newcommand{\Lie}[1][]{\operatorname{\mathcal{L}}\def\temp{#1}
\ifx\temp\empty\else^{(#1)}\fi}
\newcommand{\rmad}{\mathrm{Ad}}

\newcommand{\sympx}{{\mathfrak{X}}}
\newcommand{\sympy}{{\mathfrak{Y}}}
\newcommand{\stz}{{\mathfrak{Z}}}

\newcommand{\stka}{\mathsf{A}}
\newcommand{\stks}{\mathsf{S}}

\newcommand{\stkHom}[1][]{\mathsf{Hom}_{\raise1.5ex\hbox
to.1em{}#1}}
\newcommand{\stkMod}{\mathsf{Mod}}

\newcommand{\rmpt}{{\rm pt}}
\newcommand{\rmptt}{{\{\rm pt\}}}

\renewcommand{\to}[1][]{\xrightarrow[]{#1}}
\renewcommand{\isoto}[1][]{\xrightarrow[#1]
{{\raisebox{-.6ex}[0ex][-.6ex]{$\mspace{1mu}\sim\mspace{2mu}$}}}}

\renewcommand{\Hom}[1][]{\mathrm{Hom}_{\raise1.5ex\hbox to.1em{}#1}}
\renewcommand{\hom}[1][]{{\mathcal{H}om}_{\raise1.5ex\hbox to.1em{}#1}}
\renewcommand{\rhom}[1][]{{R\mathcal{H}om}_{\raise1.5ex\hbox to.1em{}#1}}
\renewcommand{\tens}[1][]{\mathbin{\otimes_{\raise1.5ex\hbox to-.1em{}#1}}}

\newcommand{\tw}[1]{\widetilde{#1}}

\newcommand{\for}{\mathit{for}}
\newcommand{\exten}{\mathit{ext}}

\newcommand{\eqdot}{\mathbin{:=}}

\newcommand{\cl}{\colon}
\newcommand{\scbul}{\,\raise.4ex\hbox{$\scriptscriptstyle\bullet$}\,}

\newcommand{\lp}{{\rm(}}
\newcommand{\rp}{{\rm)}}

\title{Quantization of complex Lagrangian submanifolds%
\footnotetext{\hglue-1.8em 
{\bf To appear in:} 
Adv. Math.}\footnotetext{\hglue-1.8em 
2000 AMS Mathematics Subject Classification(s): 46L65, 14A20, 32C38}
}

\author{Andrea D'Agnolo\and Pierre Schapira}

\date{}

\begin{document}

\maketitle

\begin{abstract}
Let $\Lambda$ be a smooth Lagrangian
submanifold of a complex symplectic manifold $\sympx$.
We construct twisted simple holonomic modules along $\Lambda$ in the
stack 
of  deformation-quantization modules on $\sympx$.
\end{abstract}
 
\section{Introduction}

Let $\sympy$ be a complex contact manifold. A local model for $\sympy$ is
an open subset of the projective cotangent bundle $P^*Y$ to a
complex manifold $Y$. The manifold $P^*Y$ is 
endowed with the  sheaf $\E[Y]$ of microdifferential operators
of~\cite{SKK}. 
In~\cite{K_q}, Kashiwara 
proves the existence of a canonical stack $\stkMod(\EE[\sympy])$ on
$\sympy$, 
locally equivalent to the stack of $\E[Y]$-modules. 
Let $\Lambda\subset\sympy$ be a 
smooth Lagrangian submanifold.  In the same paper, 
Kashiwara states
that there exists 
a globally defined  holonomic system 
simple along $\Lambda$ in the stack $\stkMod(\EE[\sympy]|_\Lambda)$
twisted by 
half-forms on  $\Lambda$.

Now, let $\sympx$ be a complex symplectic manifold. 
A local model for $\sympx$ is an
open subset of the  cotangent bundle $T^*X$ to a
complex manifold $X$. The manifold $T^*X$ is 
endowed with the sheaf $\W[X]$ of WKB-differential operators, 
similar to $\E[X]$, but with an extra
central parameter, a substitute to the lack of homogeneity.
(Note that,  in  the literature, $\W[X]$ is also called a 
deformation-quantization ring, or a ring of semi-classical
differential operators. See \cite{PS} for a precise
description of the ring $\W[X]$ and its
links with $\E[X]$.) A stack $\stkMod(\WW[\sympx])$ 
on $\sympx$ locally equivalent to 
the stack of $\W[X]$-modules has been constructed 
in the formal case (in the general setting of Poisson manifolds) 
by \cite{Ko} and in the analytic case (and by a different method, 
similar to~\cite{K_q})
by \cite{PS}. (See also \cite{NT0,BK,Ye} for papers closely related to
this subject.) 

In this paper, we prove that if $\Lambda$ is a smooth Lagrangian
submanifold of the complex symplectic manifold $\sympx$,
there exists a globally defined simple holonomic module
along $\Lambda$ in the stack $\stkMod(\WW[\sympx]|_\Lambda)$ twisted by
half-forms on $\Lambda$. As a by-product, we prove that there is an
equivalence of stacks between that of twisted regular holonomic
modules along $\Lambda$ and that of local systems on $\Lambda$.
The local model for our theorem is
given by $\sympx=T^*X$ and $\Lambda=T_X^*X$, the zero-section of
$T^*X$. In this case, a simple module is the sheaf $\OO_X^\tau$ whose
sections are series
$\sum_{-\infty<j \leq m}f_j\tau^j$ ($m\in\Z$), where the $f_j$'s
are sections of $\sho_X$ and the family $\{f_j\}_j$ satisfies certain
growth conditions on compact subsets of $X$. The problem we solve
here is how to patch together these local models.

Our proof consists in showing that
if $\sympx$ is a complex symplectic manifold and $\Lambda$ a  Lagrangian
submanifold,  then there exists a
``contactification'' $\sympy$ of $\sympx$ in a neighborhood of $\Lambda$.
Local models for $\sympx$ and $\sympy$ are an
open subset of the cotangent bundle $T^*X$ to a
complex manifold $X$ and an open subset of the 
projective cotangent bundle $P^*(X\times\C)$, respectively.
With the same techniques as in~\cite{K_q,PS}, we construct a
stack $\stkMod(\EE[\sympy,\hat t])$ on $\sympy$ locally equivalent
to the stack of modules over the ring $\E[X\times\C,\hat t]$ of
microdifferential operators commuting with $\partial/\partial_t$,
where $t$ is the coordinate on $\C$.
We then apply Kashiwara's existence theorem for simple
modules along Lagrangian manifolds in the contact case
to deduce the corresponding result in the symplectic case.
In fact, the technical heart of this paper is devoted to giving a
detailed proof, based on the theory of symbols of simple sections
of holonomic modules, of Kashiwara's result stated in~\cite{K_q}.

\medskip\noindent
{\bf Acknowledgements} We would like to thank 
Louis Boutet de Monvel, Masaki Kashiwara, and Pietro Polesello for their useful comments
and insights.

\section{Stacks}\label{section:stack}

Stacks were invented by Grothendieck and Giraud \cite{Gi} 
and we refer to \cite{KS2} for an exposition.
Roughly speaking, a prestack (resp.\ a stack) 
is a presheaf (resp.\ a sheaf) of categories, as we shall see below. 

In sections \ref{section:stack} and \ref{section:twmod}, 
we denote by $X$ a topological space. However, all definitions and
results easily extend when replacing $X$ with a site, that is, a small
category $\shc_X$ endowed with a Grothendieck topology.

\begin{definition}\label{def:stack1}
\begin{enumerate}[{\rm(a)}]
\item
A prestack $\stks$ on $X$ is the assignment of a category $\stks(U)$ for 
every open subset $U\subset X$, a functor 
$\rho_{VU}\cl\stks(U)\to\stks(V)$ for every open inclusion 
$V \subset U$, and an  isomorphism of functors 
$\lambda_{WVU}\cl  \rho_{WV}\rho_{VU}\Rightarrow\rho_{WU}$ 
for every open inclusion 
$W \subset V \subset U$, such that
    $\rho_{UU}=\id_{\stks(U)}$, $\lambda_{UUU} =
    \id_{\id_{\stks(U)}}$, and the following diagram of isomorphisms 
of functors from
    $\stks(U)$ to $\stks(Y)$ commutes for every open inclusion
$Y \subset W \subset V \subset U$
$$
\xymatrix 
{ \rho_{YW}  \rho_{WV}  \rho_{VU}
\ar@{=>}[rr]^-{\lambda_{YWV}  \id_{\rho_{VU}}} 
\ar@{=>}[d]^{\id_{\rho_{YW}}  \lambda_{WVU}} &&
\rho_{YV}  \rho_{VU} \ar@{=>}[d]^{\lambda_{YVU}}
\\
\rho_{YW}  \rho_{WU} \ar@{=>}[rr]^-{\lambda_{YWU}}
&& \rho_{YU}.  }
$$
For $F\in\stks(U)$ and $V\subset U$, we will write $F|_V$ for short 
instead of $\rho_{VU}(F)$.
\item
A separated prestack is a prestack $\stks$ such that for any 
$F,G\in\stks(U)$, the presheaf $\hom[\stks|_U](F,G)$, defined by 
$V\mapsto \Hom[\stks(V)](F|_V,G|_V)$, is a sheaf.
\end{enumerate}
\end{definition}

A stack is a separated prestack satisfying suitable glueing
conditions, 
which may be expressed in terms of descent data.

\begin{definition}\label{def:descdata1}
Let $U$ be an open subset of $X$, $\shu=\{U_i\}_{i\in I}$ an open
covering of $U$ and  $\stks$ a separated prestack on $X$. 
\begin{enumerate}[{\rm(a)}]
\item
A descent datum on $\shu$ for $\stks$ is a pair
\begin{equation}
\label{eq:descent}
(\{\shf_{i}\}_{i\in I}, \{\theta_{ij}\}_{i,j\in I}),
\mbox{ with } \shf_i\in\stks(U_i),
\quad \theta_{ij}\cl \shf_j\vert_{U_{ij}}\isoto\shf_i\vert_{U_{ij}}
\end{equation}
such that the following diagram
of isomorphisms in $\stks(U_{ijk})$ commutes
$$
\xymatrix{ \shf_{j}\vert_{U_{ijk}} &
\shf_{j}\vert_{U_{ij}}\vert_{U_{ijk}} \ar[l]_{\lambda}
\ar[rr]^{\theta_{ij}\vert_{U_{ijk}}} && 
\shf_{i}\vert_{U_{ij}}\vert_{U_{ijk}}
\ar[r]^{\lambda} & \shf_{i}\vert_{U_{ijk}} \\
\shf_{j}\vert_{U_{jk}}\vert_{U_{ijk}} \ar[u]^{\lambda} &&&&
\shf_{i}\vert_{U_{ik}}\vert_{U_{ijk}} \ar[u]_{\lambda} \\
& \shf_{k}\vert_{U_{jk}}\vert_{U_{ijk}} \ar[r]^-{\lambda}
\ar[ul]^{\theta_{jk}\vert_{U_{ijk}}} & \shf_{k}\vert_{U_{ijk}} &
\shf_{k}\vert_{U_{ik}}\vert_{U_{ijk}}. 
\ar[ur]_{\theta_{ik}\vert_{U_{ijk}}}
\ar[l]_-{\lambda} & } $$
\item
The descent datum \eqref{eq:descent} is called
effective if  there exist $\shf \in\stks(U)$ 
and isomorphisms 
$\theta_i\cl \shf\vert_{U_i}\isoto\shf_i$ in $\stks(U_i)$
satisfying the natural compatibility conditions with the
$\theta_{ij}$'s and $\lambda$'s.
\end{enumerate}
\end{definition}
Note that if the descent datum \eqref{eq:descent} is effective, 
then $\shf$ is  unique up to 
unique isomorphism.

\begin{definition}
A stack is a separated prestack such that for any 
 open subset $U$ of $X$ and any open covering 
$\shu=\{U_i\}_{i\in I}$ of $U$, the 
descent datum  is effective.
\end{definition}
To end this section, let us go up one level, and recall the glueing 
conditions for stacks.

\begin{definition}\label{def:descdata2}
Let $U$ be an open subset of $X$, $\shu=\{U_i\}_{i\in I}$ an open
covering of $U$. 
\begin{enumerate}[{\rm(a)}]
\item
A descent datum for stacks on  $\shu$
is a triplet
\begin{equation}
\label{eq:stksi}
( \{\stks_i\}_{i\in I}, \{\varphi_{ij}\}_{i,j\in I}, 
\{\alpha_{ijk}\}_{i,j,k\in I} ),
\end{equation}
where the 
$\stks_i$'s are stacks on $U_i$,
$\varphi_{ij}\cl\stks_j\vert_{U_{ij}} 
\to \stks_i\vert_{U_{ij}}$ are equivalences of stacks, and
$\alpha_{ijk}\cl \varphi_{ij} \varphi_{jk}\Rightarrow
\varphi_{ik} \cl \stks_k\vert_{U_{ijk}} \to \stks_i\vert_{U_{ijk}}$ 
are isomorphisms of functors such that for any $i,j,k,l \in I$ 
the following diagram of isomorphisms of functors from 
$\stks_l|_{U_{ijkl}}$ to $\stks_i|_{U_{ijkl}}$ commutes
\begin{equation}
\label{eq:phiijcommut}
\xymatrix@C=4em{ {\varphi_{ij}\varphi_{jk}\varphi_{kl}}
\ar@{=>}[r]^-{\alpha_{ijk}  \id_{\varphi_{kl}}} 
\ar@{=>}[d]^{\id_{\varphi_{ij}}  \alpha_{jkl}}
&{\varphi_{ik}\varphi_{kl}}\ar@{=>}[d]^{\alpha_{ikl}}\\
{\varphi_{ij}\varphi_{jl}}\ar@{=>}[r]^{\alpha_{ijl}}&\varphi_{il}.
}
\end{equation}
\item
The descent datum \eqref{eq:stksi} is called effective if there
exist a stack $\stks$ on $U$, equivalences 
of stacks $\varphi_i\cl \stks\vert_{U_i}\to\stks_i$ and 
isomorphisms of functors 
$\alpha_{ij}\cl \varphi_{ij}\varphi_j
\Rightarrow\varphi_i\cl\stks\vert_{U_{ij}} 
\to \stks_i\vert_{U_{ij}}$,  
satisfying the natural compatibility conditions. 
\end{enumerate}
\end{definition}
Note that if the descent datum \eqref{eq:phiijcommut} is effective, 
then $\stks$  is  unique up to 
equivalence and such an equivalence is unique up to unique
isomorphism.

In the language of $2$-categories, the following theorem asserts that
{\em the $2$-prestack of stacks is a $2$-stack}.

\begin{theorem}\label{th:desstk}
{\rm (cf~\cite{Gi,Br})} Descent data for stacks are effective.
\end{theorem}

Denote by $\stks$ the stack associated with
the descent datum \eqref{eq:phiijcommut}. Objects of $\stks(U)$ can
be 
described by pairs
\begin{equation}
\label{eq:shfi}
( \{\shf_i\}_{i\in I}, \{\xi_{ij}\}_{i,j\in I} ),
\end{equation}
where $\shf_i\in \stks_i(U_i)$ and
$\xi_{ij}\cl\varphi_{ij}(\shf_{j}\vert_{U_{ij}})\to\shf_i\vert_{U_{ij}}$
are isomorphisms in $\stks_i(U_{ij})$ such that
for $i,j,k\in I$ the following diagram in $\stks_i(U_{ijk})$ commutes
\begin{equation}
\label{eq:xixiaxi}
\xymatrix{ 
\varphi_{ij}  \varphi_{jk} ( \shf_{k}\vert_{U_{ijk}})
\ar[r]^-{\varphi_{ij}(\xi_{jk})} \ar[d]^{\alpha_{ijk}(\shf_k)} &
\varphi_{ij}(\shf_{j}\vert_{U_{ijk}}) 
\ar[d]^{\xi_{ij}} \\
\varphi_{ik}(\shf_{k}\vert_{U_{ijk}}) 
\ar[r]^-{\xi_{ik}} &
\shf_{i}\vert_{U_{ijk}}.}
\end{equation}

\section{Twisted modules}\label{section:twmod}

Let us now recall how stacks of twisted modules are constructed 
in~\cite{K_tw,K_q}. 

Let $\field$ be a commutative unital ring and $\sha$ a 
sheaf of $\field$-algebras on $X$. 
Denote by $\catMod(\sha)$ the category of left $\sha$-modules, 
and by $\stkMod(\sha)$ the associated stack 
$X\supset U \mapsto \catMod(\sha|_U)$.

Consider an open covering $\shu=\{U_i\}_{i\in I}$ of $X$, a family 
of $\field$-algebras $\sha_i$ on $U_i$ and $\field$-algebra 
isomorphisms 
$f_{ij}\cl \sha_j|_{U_{ij}} \to\sha_i|_{U_{ij}}$. 
The existence of a sheaf of $\field$-algebras
locally isomorphic to $\sha_i$ requires the condition $f_{ij}f_{jk} =
f_{ik}$ on triple intersections. The weaker conditions
\eqref{eq:ffadf} 
and \eqref{eq:aafaa} below are needed for the existence of a 
$\field$-additive stack locally equivalent to $\stkMod(\sha_i)$. 

\begin{definition}\label{de:alg}
A $\field$-algebroid descent datum $\stka$ on $\shu$ 
is a triplet
\begin{equation}
\label{eq:shai}
\stka = (\{\sha_i\}_{i\in I},\{f_{ij}\}_{i,j \in
I},\{a_{ijk}\}_{i,j,k \in I}),
\end{equation}
where $\sha_i$ is a  $\field$-algebra on $U_i$, 
$f_{ij}\cl \sha_j|_{U_{ij}}\to\sha_i|_{U_{ij}}$ 
is a $\field$-algebra isomorphism, 
$a_{ijk} \in \sha_i^\times(U_{ijk})$ is an invertible section, 
and \eqref{eq:ffadf} and \eqref{eq:aafaa} below are satisfied:
\begin{align}
\label{eq:ffadf}
f_{ij}f_{jk} &= \rmad(a_{ijk}) f_{ik} \text{ as $\field$-algebra 
isomorphisms $\sha_k|_{U_{ijk}} \isoto \sha_i|_{U_{ijk}}$}, \\
\label{eq:aafaa}
a_{ijk} a_{ikl} &= f_{ij}(a_{jkl}) a_{ijl}
\mbox{ in }\sha_i^\times(U_{ijkl}).
\end{align}
(Here $\rmad(a_{ijk})$ denotes the automorphism of $\sha_i|_{U_{ijk}}$ given by $a\mapsto a_{ijk}\, a \,\opb{a_{ijk}}$.)
\end{definition}

\begin{remark}
The notion of an algebroid stack exists intrinsically,
without using coverings or descent data. It has been introduced 
by~\cite{Ko} and developed in~\cite{DP2}. 
In this paper, we shall restrict ourselves to algebroids presented by 
descent data.
\end{remark}

\begin{remark}
\label{re:aafaa}
Let $\shb_i\subset\sha_i^\times$ be multiplicative subgroups,
invariant 
by $f_{ij}$, and such that for any $b_i,b_i'\in\shb_i$, the equality
$\mathrm{Ad}(b_i) = \mathrm{Ad}(b_i')$ implies $b_i = b_i'$. 
Assume that $a_{ijk}\in\shb_i$. Then, as noticed {\em e.g.},
 in \cite[pag.~2]{K_q}, 
condition \eqref{eq:aafaa} follows from \eqref{eq:ffadf}.
\end{remark}
Let us recall how to define the stack of ``$\stka$-modules'' 
 in terms of local data.

A $\field$-algebra morphism $f\cl\shb\to\sha$ induces a functor 
$$
\tw f \cl \stkMod(\sha) \to \stkMod(\shb)
$$ 
defined by $\shm \mapsto \atw f \shm$, where $\atw f \shm$ denotes the
sheaf of $\field$-vector spaces $\shm$, endowed with the $\shb$-module structure 
given by $b m\eqdot f(b) m$ for $b\in \shb$ and $m \in \shm$.

For $a\in\sha^\times$ an invertible section, the automorphism 
$\rmad(a)$ induces the functor $\widetilde{\rmad(a)}$
between $\stkMod(\sha)$ and itself, and we denote by 
$$
\tw a \cl \widetilde{\rmad(a)} \Rightarrow \id_{\stkMod(\sha)}
$$ 
the isomorphism of functors given by 
$\tw a(\shm)\cl  \atw{\rmad(a)}\shm \to \shm$,  
$u\mapsto \opb{a}u$, for $u\in\shm\in\stkMod(\sha)$. (Note that 
$\tw a(\shm)(a'u)=\opb{a}aa'\opb{a}u=a'\tw a(\shm)(u)$.)

\begin{definition}\label{de:algmod}
\begin{itemize}
\item[(i)]
The stack of twisted modules associated to the $\field$-algebroid
descent datum $\stka$ on $\shu$ in \eqref{eq:shai} 
is the stack defined (using Theorem \ref{th:desstk})
by the descent datum
\begin{equation}
\label{eq:stkmoda}
\stkMod(\stka) = (\{\stkMod(\sha_i)\}_{i\in I},\{\tilde f_{ji}\}_{i,j\in
I},
\{\tilde a_{kji}\}_{i,j,k \in I}).
\end{equation}
\item[(ii)]
One sets $\catMod(\stka) \eqdot \stkMod(\stka)(X)$. Objects of the
category $\catMod(\stka)$
are called twisted modules.
\end{itemize}
\end{definition}

According to \eqref{eq:shfi}, objects of $\catMod(\stka)$ 
are described by pairs
$$
\shm = ( \{\shm_i\}_{i \in I}, \{\xi_{ij}\}_{i,j \in I} ),
$$
where $\shm_i$ are $\sha_i$-modules and
$\xi_{ij}\cl \atw{f_{ji}}\shm_j\vert_{U_{ij}} \to \shm_{i}
\vert_{U_{ij}}$ are isomorphisms of $\sha_i$-modules, such that 
for any $u_k\in\shm_k$ one has
\begin{equation}
\label{eq:glueM}
\xi_{ij}  ({}_{f_{ji}}\xi_{jk} (u_k))
= \xi_{ik}(a_{kji}^{-1}  u_k)
\end{equation}
as morphisms $\atw{f_{kj}f_{ji}}\shm_k \to \shm_i$.
Indeed, \eqref{eq:glueM} translates the commutativity of
\eqref{eq:xixiaxi}.

\begin{example}\label{ex:Omega1/2}
Let $X$ be a complex manifold, and denote by $\Omega_X$ 
the sheaf of holomorphic forms of maximal degree. 
Take an open covering 
$\shu=\{U_i\}_{i\in I}$ of $X$ such that there are nowhere vanishing 
sections $\omega_i\in\Omega_{U_i}$. 
Let $t_{ij} \in \OO^\times_{U_{ij}}$ be the 
transition functions given by 
$\omega_j|_{U_{ij}} = t_{ij} \omega_i|_{U_{ij}}$. 
Choose determinations $s_{ij} \in \OO^\times_{U_{ij}}$ 
for the multivalued functions $t_{ij}^{1/2}$. 
Since $s_{ij}s_{jk}$ and $s_{ik}$ are both determinations 
of $t_{ik}^{1/2}$, there exists $c_{ijk} \in\{-1,1\}$ such that
\begin{equation}
\label{eq:sscs}
s_{ij}s_{jk} = c_{ijk}s_{ik}.
\end{equation}
We thus get a $\C$-algebroid descent datum
$$
\C_{\sqrt{\Omega_X}}\eqdot ( \{ \C_{U_i} \}_{i \in I}, 
\{\id_{\C_{U_{ij}} } \}_{i,j \in I}, \{c_{ijk}\}_{i,j,k \in I}).
$$
Note that, since $c_{ijk}^2=1$, there is an equivalence 
$\catMod(\C_{\sqrt{\Omega_X}}) \simeq \catMod(\C_{\sqrt{\Omega_X^{-1}}})$.

Recall from \cite{K_tw} (see also~\cite[\S1]{DS}), that  there is an equivalence 
\begin{equation}
\label{eq:trivtw}
\catMod(\C_{\sqrt{\Omega_X}}) \simeq \catMod(\C_X)
\end{equation}
if and only if the cohomology class $[c_{ijk}]\in H^2(X;\C_X^\times)$ is trivial.
Consider the long exact cohomology sequence
$$
H^1(X;\C_X^\times) \to[\alpha]
H^1(X;\OO_X^\times) \to[\beta]
H^1(X;d\OO_X) \to[\gamma]
H^2(X;\C_X^\times)
$$
associated with the short exact sequence
$$
1 \to \C_X^\times \to \OO_X^\times \to[d\log] d\OO_X \to 0.
$$
One has $[c_{ijk}] = \gamma(\frac12\beta([\Omega_X]))$, so that $[c_{ijk}] =1$ if and only if 
there exists a line bundle $\shl$ such that $\frac12\beta([\Omega_X]) = \beta([\shl])$, i.e.~such that
$\beta([\Omega_X\tens[\OO]\shl^{-2}]) = 0$. This last condition holds if and only if there exists a local 
system of rank one $L$ such that $[\Omega_X\tens[\OO]\shl^{-2}] = \alpha([L])$. Summarizing,
\eqref{eq:trivtw} holds if and only if there exist $\shl$ and $L$ as above, such that
$$
\Omega_X \simeq L \tens \shl^2.
$$

The twisted sheaf of half-forms in $\catMod(\C_{\sqrt{\Omega_X}})$ is given by
$$
\sqrt{\Omega_X} = ( \{ \OO_{U_i} \}_{i \in I}, 
\{ s_{ij} \}_{i,j \in I} ). 
$$
We denote by $\sqrt{\omega_i}$ the section corresponding to $1\in \OO_{U_i}$. 
Hence, on $U_{ij}$ we have
$$
\sqrt{\omega_j} = s_{ij} \sqrt{\omega_i}.
$$
Denote by $[s_{ij}]$ the equivalence class of $s_{ij}$ in
$\OO^\times_{U_{ij}}/\C_{U_{ij}}^\times$. Since the $s_{ij}$'s
satisfy \eqref{eq:sscs}, we notice that
\begin{equation}
\label{eq:Omega12C}
\sqrt{\Omega_X^\times}/\C_X^\times = ( \{
\OO^\times_{U_i}/\C_{U_i}^\times \}_{i \in I}, 
\{ [s_{ij}] \}_{i,j \in I} ) \in \catMod(\Z_X)
\end{equation}
is a (usual, i.e.~not twisted) sheaf.
\end{example}

Let 
$\stka=(\{\sha_i\}_{i\in I},\{f_{ij}\}_{i,j\in I},\{a_{ijk}\}_{i,j,k
\in I})$ 
be a $\field$-algebroid descent datum on $\shu$ as in
\eqref{eq:shai}, 
and consider a pair
$$
\shm = ( \{\shm_i\}_{i \in I}, \{\xi_{ij}\}_{i,j \in I} ),
$$
where $\shm_i$ are $\sha_i$-modules and
$\xi_{ij}\cl \atw{f_{ji}}\shm_j\vert_{U_{ij}} \to \shm_{i}
\vert_{U_{ij}}$ are isomorphisms of $\sha_i$-modules which
do not necessarily satisfy \eqref{eq:glueM}. Assume instead that
there are isomorphisms
\begin{equation}
\label{eq:homK}
\hom[\sha_i](\atw{f_{ji}}\shm_j\vert_{U_{ij}},\shm_{i}\vert_{U_{ij}})
\simeq \field_{U_{ij}}.
\end{equation}
Under this assumption, we will show in
Proposition~\ref{pr:twistedglue} below that $\shm$ 
makes sense as a global object of $\stkMod(\stka\tens[\field]\stks)$, for a
suitable twist $\stks$.

Consider the isomorphisms
\begin{equation}
\label{eq:phiji}
\phi_{ij}\colon
\hom[\sha_i](\atw{f_{ji}}\shm_j\vert_{U_{ij}},\shm_{i}\vert_{U_{ij}})
\isoto \field_{U_{ij}},
\end{equation}
defined by $\xi_{ij} \mapsto 1$.
The dotted arrow defined by the following commutative diagram of
isomorphisms is the multiplication by a section $c_{ijk}$ of
$\field_{U_{ijk}}^\times$. 
\begin{equation}
\label{eq:cijk}
\xymatrix{
\hom[\sha_j](\atw{f_{ij}}\shm_i, \shm_j) \tens[\field]
\hom[\sha_k](\atw{f_{jk}}\shm_j, \shm_k) 
\ar[r]_-\sim^-{\varphi_{ji}\tens\varphi_{kj}}
\ar[d]_\wr^{\tw{f_{jk}}\tens\id} &
\field_{U_{ijk}}\tens[\field]\field_{U_{ijk}} \ar[d]_\wr \\
\hom[\sha_k](\atw{f_{ij}f_{jk}}\shm_i, \atw{f_{jk}}\shm_j)
\tens[\field]
\hom[\sha_k](\atw{f_{jk}}\shm_j, \shm_k) & 
\field_{U_{ijk}} \\
\hom[\sha_k](\atw{f_{ik}}\shm_i, \atw{f_{jk}}\shm_j) \tens[\field]
\hom[\sha_k](\atw{f_{jk}}\shm_j, \shm_k) 
\ar[d]_\wr^\circ  \ar[u]^\wr_{\tw{a_{ijk}}\tens\id} \\
\hom[\sha_k](\atw{f_{ik}}\shm_i, \shm_k) \ar[r]_-\sim^-{\varphi_{ki}}
&
\field_{U_{ijk}}. \ar@{.>}[uu]_\wr^{c_{ijk}}
}
\end{equation}
Here, the first vertical arrow on the left follows from the equality
$\atw{f_{ij}f_{jk}}\shm_i = \atw{f_{jk}}(\atw{f_{ij}}\shm_i)$, and
the second one from $\atw{f_{ij}f_{jk}}\shm_i =
\atw{\rmad(a_{ijk})f_{ik}}\shm_i$.

\begin{lemma}\label{le:cocycle1}
The constants $c_{ijk}$ defined above satisfy the cocycle condition
$$
c_{ijk} c_{ikl} = c_{jkl} c_{ijl}.
$$
\end{lemma}

\begin{proof}
Consider morphisms
$$
\eta_{ji} \colon \atw{f_{ij}}\shm_i \to \shm_j.
$$
These induce morphisms
$$
\atw{f_{jk}}\eta_{ji} \colon \atw{f_{jk}}(\atw{f_{ij}}\shm_i)  \to
\atw{f_{jk}}\shm_j.
$$
The composition
$$
\atw{f_{ik}}\shm_i \to[\tw{a_{ijk}}^{-1}]
\atw{\rmad(a_{ijk})f_{ik}}\shm_i = \atw{f_{jk}}(\atw{f_{ij}}\shm_i)
\to[\atw{f_{jk}}\eta_{ji}] \atw{f_{jk}}\shm_j
$$
is given by
$$
u_i \mapsto a_{ijk} u_i \mapsto \atw{f_{jk}}\eta_{ji}(a_{ijk}u_i).
$$
Hence, the composition of the vertical isomorphisms in the left 
column of \eqref{eq:cijk} is given by
$$
\eta_{ji}\tens\eta_{kj} \mapsto \eta_{kj}
(\atw{f_{jk}}\eta_{ji}(a_{ijk}\cdot\scbul)).
$$ 
One then has
\begin{equation}
\label{eq:phijikjki}
\varphi_{ji}(\eta_{ji})\varphi_{kj}(\eta_{kj})=
\varphi_{ki}(\eta_{kj} (\atw{f_{jk}}\eta_{ji}(a_{ijk}\cdot\scbul)))
c_{ijk}.
\end{equation}
Using \eqref{eq:phijikjki}, we have, on one hand
\begin{align*}
\varphi_{ji}(\eta_{ji})\varphi_{kj}(\eta_{kj})\varphi_{lk}(\eta_{lk}) 
&=\varphi_{ji}(\eta_{ji})\varphi_{lj}(\eta_{lk}(\atw{f_{kl}}\eta_{kj}
(a_{jkl}\cdot\scbul)))
c_{jkl} \\
&=\varphi_{li}(\eta_{lk}(\atw{f_{kl}}\eta_{kj}(a_{jkl}\atw{f_{jl}}
\eta_{ji}(a_{ijl}\cdot\scbul))))c_{ijl} 
c_{jkl} \\
&=\varphi_{li}(\eta_{lk}(\atw{f_{kl}}\eta_{kj}(\atw{f_{jk}f_{kl}}
\eta_{ji}(f_{ij}(a_{jkl})a_{ijl}\cdot\scbul))))
c_{ijl} c_{jkl},
\end{align*}
and on the other hand
\begin{align*}
\varphi_{ji}(\eta_{ji})\varphi_{kj}(\eta_{kj})\varphi_{lk}(\eta_{lk})
&=
\varphi_{ki}(\eta_{kj}(\atw{f_{jk}}\eta_{ji}(a_{ijk}\cdot\scbul)))
\varphi_{lk}(\eta_{lk}) c_{ijk} \\
&= \varphi_{li}(\eta_{lk}(\atw{f_{kl}}\eta_{kj}(\atw{f_{jk}f_{kl}}
\eta_{ji}(a_{ikl}a_{ijk}\cdot\scbul)))) c_{ikl} c_{ijk}.
\end{align*}
The conclusion follows using \eqref{eq:aafaa}.
\end{proof}

\begin{remark}
Lemma \ref{le:cocycle1} is a particular case of a general result which
asserts that equivalence classes of locally trivial
$\field$-algebroids on
$X$ are
in one-to-one correspondence with $H^2(X;\field^\times)$. 
The analogue result for gerbes is discussed in~\cite{Gi}, 
and we refer to \cite{DP2} for the formulation in terms of algebroids.
\end{remark}

Let us recall our setting. Given a $\field$-algebroid descent datum
$$
\stka = (\{\sha_i\}_{i\in I},\{f_{ij}\}_{i,j \in
I},\{a_{ijk}\}_{i,j,k \in I}),
$$
consider a pair
$$
\shm = ( \{\shm_i\}_{i \in I}, \{\xi_{ij}\}_{i,j \in I} ),
$$
where $\shm_i$ are $\sha_i$-modules and
$\xi_{ij}\cl \atw{f_{ji}}\shm_j\vert_{U_{ij}} \to \shm_{i}
\vert_{U_{ij}}$ are isomorphisms of $\sha_i$-modules. Assuming
\eqref{eq:homK}, Lemma~\ref{le:cocycle1} guarantees that
the constants $c_{ijk}$ defined by \eqref{eq:cijk} satisfy the
cocycle condition.
We can thus consider the $\field$-algebroid descent datum
$$
\stks = ( \{ \field_{U_i} \}_{i \in I}, 
\{\id_{\field_{U_{ij}} } \}_{i,j \in I}, \{c_{ijk}\}_{i,j,k \in I}).
$$
Set
$$
\stka\tens[\field]\stks = 
(\{\sha_i\}_{i\in I},\{f_{ij}\}_{i,j \in I},\{a_{ijk}c_{ijk}\}_{i,j,k
\in I}).
$$

\begin{proposition}
\label{pr:twistedglue}
Let $\shm = ( \{\shm_i\}_{i \in I}, \{\xi_{ij}\}_{i,j \in I} )$ be as
above, and assume \eqref{eq:homK}, then
$$
\shm = (\shm_i,\xi_{ij})\in \catMod(\stka\tens[\field]\stks).
$$
\end{proposition}

\begin{proof}
By \eqref{eq:glueM}, it is enough to show that
$$
\xi_{ik} = \xi_{ij}(\atw{f_{ji}}\xi_{jk}(a_{kji}c_{kji}\cdot\scbul)).
$$
We have
$\varphi_{ik}(\xi_{ik}) = 1 =
\varphi_{jk}(\xi_{jk})\varphi_{ij}(\xi_{ij})=
\varphi_{ik}(\xi_{ij}(\atw{f_{ji}}\xi_{jk}(a_{kji}c_{kji}\cdot\scbul)))$
by
\eqref{eq:phijikjki}. 
The conclusion follows since  $\varphi_{ik}$ is an isomorphism.
\end{proof}

\section{Microdifferential modules}\label{se:mumo}

Here we review a few notions from the theory of 
microdifferential modules. References are made to \cite{SKK} and
also to  \cite{KO} for complementary results. 
See  \cite{K3, Sc} for an exposition.

Let $Y$ be a complex analytic manifold, and $\pi\cl T^*Y \to Y$ 
its cotangent bundle.
The sheaf $\E[Y]$  of microdifferential operators on $T^*Y$ is a 
$\C$-central algebra endowed with a $\Z$-filtration by the order.
Denote by $\E[Y](m)$ its subsheaf of operators of order at 
most $m$, and by $\OO_{T^*Y}(m)$ the sheaf of functions 
homogeneous of degree $m$ in the fiber of $\pi$. Denote by
$\operatorname{eu}$ the Euler vector field, i.e.~the infinitesimal
generator of the action of $\C^\times$ on $T^*Y$. Then $f\in
\OO_{T^*Y}(m)$ if and only if $\operatorname{eu} f = m f$.
In a homogeneous symplectic 
local coordinate system  $(x;\xi)$ on
$U\subset T^*Y$, a section $P\in \sect(U;\E[Y](m))$ 
is written as a formal series
\begin{equation}
\label{eq:P}
P = \sum_{j\leq m}p_j(x;\xi), \quad p_j\in\sect(U;\OO_{T^*Y}(j)),
\end{equation}
with the condition that for any compact subset $K$ of $U$ there
exists a 
constant $C_K>0$ such that 
$\sup\limits_{K}\vert p_{j}\vert \leq C_K^{-j}(-j)!$ for all $j<0$.

If $Q=\sum_{j\leq n}q_j$ is another section, 
the product $PQ = R = \sum_{j\leq m+n}r_j$ is given by the Leibniz
rule
$$
r_k = \sum_{k=i+j-|\alpha|} \frac 1{\alpha!} 
(\partial_\xi^\alpha p_i) (\partial_x^\alpha q_j).
$$

The symbol map
$$
\sigma_m \cl \E[Y](m) \to \OO_{T^*Y}(m), \qquad P\mapsto p_m
$$
does not depend on the choice of coordinates and 
induces the symbol map 
$$
\sigma\cl\E[Y]\to \gr\E[Y] \isoto \DSum_{m \in \Z}\OO_{T^*Y}(m).
$$
The formal adjoint of $P = \sum_{j\leq m}p_j$ is defined by
$$
P^* = \sum_{j\leq m}p^*_j, \quad p^*_j(x;\xi) 
=\sum_{j=k-|\alpha|}\frac{(-1)^{|\alpha|}}{\alpha!}
\partial_\xi^\alpha \partial_x^\alpha p_k(x;-\xi).
$$
It depends on the choice of coordinates, and more precisely on 
the choice of the top degree form 
$dx_i\wedge\cdots \wedge dx_n\in\Omega_Y$. One thus 
considers the twist of $\E[Y]$ by half-forms
$$
\Ev[Y] \eqdot \opb\pi\sqrt{\Omega_Y} \tens[\opb\pi\OO_Y] \E[Y] 
\tens[\opb\pi\OO_Y] \opb\pi\sqrt{\Omega_Y^{-1}}.
$$ 
This is a sheaf of filtered $\C$-algebras endowed with a canonical
anti-isomorphism
$$
*\cl \Ev[Y] \isoto \oim a\Ev[Y],
$$
where $a$ denotes the antipodal map on $T^*Y$.
There is a  subprincipal symbol
$$
\sigma'_{m-1} \cl \Ev[Y](m) \to \OO_{T^*Y}(m-1), \quad
P \mapsto \frac12\sigma_{m-1}(P- (-1)^m P^*).
$$
In local coordinates, 
$\sigma'_{m-1}(P)
=p_{m-1}-\frac12\sum\nolimits_i\partial_{x_i}\partial_{\xi_i}p_m$ 
for $P$ as in \eqref{eq:P}. 

The following definition is adapted from  \cite{KO}.

\begin{definition}\label{def:simpleEmod}
Let $\Lambda$ be a smooth, locally closed, 
$\C^\times$-conic submanifold of $T^*Y$. 
Let $\shm$ be a coherent $\E[Y]$-module supported by $\Lambda$.  
\begin{enumerate}[{\rm(a)}]
\item
One says that $\shm$ is regular (resp.~simple) along
$\Lambda$ if there locally exists a coherent sub-$\E[Y] (0)$-module
$\shm_0$
of $\shm$ which generates it over $\E[Y] $, and such that
$\shm_0/\E[Y] (-1)\shm_0$ is an $\OO_\Lambda(0)$-module (resp.~a
locally free $\OO_\Lambda(0)$-module of
rank one).  
\item
Let $\shm$ be simple along $\Lambda$. A
section $u\in\shm$ is called a simple generator if 
$\shm_0=\E[Y](0)u$ satisfies the conditions in (a), and the image of
$u$ in $\shm_0/\E[Y](-1)\shm_0$
generates this module over $\OO_\Lambda(0)$.
\end{enumerate}
\end{definition}

Set
\begin{equation}
\label{eq:shilambda}
\shi_{\Lambda} = \{ P\in\E[Y](1)|_\Lambda ;\ \sigma_1(P)|_\Lambda = 0
\},
\end{equation}
and denote by $\E[\Lambda]$ the sub-algebra of $\E[Y]|_\Lambda$ 
generated by $\shi_{\Lambda}$.

\begin{remark}
Let $u$ be a generator of a coherent $\E[Y]$-module $\shm$. 
Then $\shm\simeq \E[Y]/\shi$, where 
$\shi = \{P\in\E[Y] \cl Pu=0\}$. Set 
$$
\overline\shi=\{\sigma_m(P); m\in\Z,\ P\in\shi\cap\E[Y](m)\}, 
$$
and note that $\supp\shm = \supp(\OO_{T^*Y}/\overline\shi)$. 
Then $\shm$ is simple if and only if 
there locally exists a generator $u$ such that the ideal $\overline
\shi$ is reduced. Moreover, such a section $u$ is a simple
generator
and  the sub-$\E[Y](0)$-module $\shm_0$ generated by $u$ satisfies
$$
\E[\Lambda]\shm_0\subset\shm_0.
$$
Indeed, in a homogeneous symplectic local coordinate system we may write $P\in\shi_{\Lambda}$ as $P=P'+Q$ with $Q$ of 
order $\leq 0$ and $P'u=0$. 
\end{remark}

\subsection*{Symbol of sections of simple systems}

Let us recall the notion of symbol for simple generators. References
are made to~\cite{SKK,K_mu}. 

For a vector field $v\in\Theta_\Lambda$  on $\Lambda$, 
denote by $\Lie^{1/2}_v$ its Lie derivative action on the 
twisted sheaf $\sqrt{\Omega_\Lambda}$. Then $\Lie^{1/2}_v$ 
is an operator of order one in the ring 
$\D[\Lambda]^\halfform = 
\sqrt{\Omega_\Lambda} 
\tens[\OO_\Lambda]\D[\Lambda]\tens[\OO_\Lambda]\sqrt{\Omega_\Lambda^{-1}}$

of differential operators acting on $\sqrt{\Omega_\Lambda}$. 

Define $\Iv[\Lambda]$ and $\Ev[\Lambda]$ as in \eqref{eq:shilambda},
replacing $\E[Y]$ with $\Ev[Y]$, and denote by  $H_f$ the 
Hamiltonian vector field of $f\in\OO_{T^*Y}$. Note that $H_f\in
T\Lambda$ if $f|_\Lambda=0$.
For $P\in \Iv[\Lambda]$, consider the transport operator
$$
\Lrm(P)=\Lie^{1/2}_{H_{\sigma_1(P)}|_\Lambda}+\sigma'_0(P)|_\Lambda.
$$
One checks that $\Lrm$ satisfies the relations
$\Lrm(AP) = \sigma_0(A) \Lrm(P)$, 
$\Lrm(PA) = \Lrm(P) \sigma_0(A)$, and
$\Lrm([P,Q]) = [ \Lrm(P) , \Lrm (Q) ]$, for $P,Q\in\Iv[\Lambda]$ and
$A \in \Ev[\Lambda](0)$ (see {\em e.g.}~\cite[\S8.3]{K3}).
It follows that $\Lrm$ extends as a $\C$-algebra morphism
\begin{equation}
\label{eq:preL}
\Lrm \cl \Ev[\Lambda] \to \D[\Lambda]^\halfform
\end{equation}
by setting $\Lrm (P_1\cdots P_r) = \Lrm (P_1) \cdots \Lrm (P_r)$, 
for $P_i \in \Iv[\Lambda]$. 

Let $\shm$ be a simple $\Ev[Y]$-module along $\Lambda$, and
$u\in\shm$ 
a simple generator. The twisted subsheaf of $\sqrt{\Omega_\Lambda}$
defined by
\begin{equation}
\label{eq:defsymb}
\{\sigma\in\sqrt{\Omega_\Lambda}
;\ \Lrm(P)\sigma = 0\ \forall P\in\Ev[\Lambda],\ Pu=0\}
\end{equation}
is locally a free sheaf of rank one over $\C$.

\begin{definition}
Let $\shm$ be a simple $\Ev[Y]$-module along $\Lambda$ and let
$u\in\shm$ 
be a simple generator. The symbol 
of $u$ is defined by
$$
\sigma_\Lambda(u) = [\sigma] \in
\sqrt{\Omega_\Lambda^\times}/\C_\Lambda^\times.
$$
for $\sigma$ as in \eqref{eq:defsymb}.
\end{definition}

The Euler vector field $\operatorname{eu}$ acts on
$\sqrt{\Omega_\Lambda}$ by $\Lie^{1/2}_{\operatorname{eu}}$, and one
says that $\sigma\in\sqrt{\Omega_\Lambda}$ is homogeneous of degree
$\lambda$ if $\operatorname{eu}\sigma = \lambda\sigma$. Hence, the
notion of homogeneous section of
$\sqrt{\Omega_\Lambda^\times}/\C_\Lambda^\times$ makes sense. One
calls
{\em order of $u$} the homogeneous degree of $\sigma_\Lambda(u)$.
Then 
the equivalence class of $\lambda$ in $\C/\Z$ does not depend on $u$,
and is 
called the {\em order of $\shm$}.

If $P\in \Ev[Y](m)$ is such that $\sigma_m(P)|_\Lambda$ 
never vanishes, then
\begin{equation}
\label{eq:sigmaPu}
\sigma_\Lambda(P u) = \sigma_m(P) \sigma_\Lambda(u).
\end{equation}
Also recall from {\em loc.cit.} that simple modules of the same order are locally
isomorphic.

\section{Quantization of contact manifolds}

Here we review Kashiwara's construction \cite{K_q} of 
the stack of microdifferential modules on a contact manifold.

Let $Y$ be a complex analytic manifold, $\pi\cl T^*Y \to Y$ 
its cotangent bundle, $\dot T^*Y = T^*Y\setminus Y$ the complementary 
of the zero-section, $\varpi\cl P^*Y \to Y$ the projective
cotangent 
bundle, and $\gamma \cl \dot T^*Y \to P^*Y$ the projection. 
The sheaf of microdifferential operators on $P^*Y$ is given by 
$\oim\gamma(\E[Y]|_{\dot T^*X})$, and we still denote it 
by $\E[Y]$ for short. Since the antipodal map induces the 
identity on $P^*Y$, the anti-involution $*$ is  well 
defined on the sheaf $\Ev[Y]$. 

Let $\chi_{ij}\cl P^*Y_i\supset V_i \to V_j\subset P^*Y_j$ 
be a contact transformation.
Recall that there locally exists a quantized contact transformation 
(QCT for short) above $\chi_{ij}$. This is an isomorphism of filtered 
$\C$-algebras 
$\Phi_{ij}\cl \chi_{ij}^{-1}\Ev[Y_j]|_{V_j} \isoto \Ev[Y_i]|_{V_i}$. 
Moreover, one can ask that $\Phi_{ij}$ is $*$-preserving. Such a 
quantization is not unique, but a key remark by Kashiwara 
is that $*$-preserving filtered automorphisms of $\Ev[Y]$ 
are of the form $\mathrm{Ad}(Q)$ for a {\em unique} 
operator $Q\in\Ev[Y](0)$ satisfying
\begin{equation}
\label{eq:Qcond}
QQ^* = 1, \quad \sigma_0(Q)=1.
\end{equation}

\begin{definition}
\label{def:contactmnf}
%(see~\cite{K_q})
A complex contact manifold $\sympy = (\sympy,\sho_\sympy(1),\contatto)$ is
a complex manifold $\sympy$ of dimension $2n+1$ endowed 
with a line bundle $\sho_\sympy(1)$ and a 1-form 
$\contatto\in\sect(\sympy,\Omega^1_\sympy\tens[\sho]\sho_\sympy(1))$, such 
that $\contatto\wedge(d\contatto)^{\wedge n}$ is a non-degenerate 
section of $\Omega^{2n+1}_\sympy\tens[\sho]\sho_\sympy(n+1)$. 
\end{definition}
(Here we set $\sho_\sympy(k) = \sho_\sympy(1)^{\otimes k}$, and we use
the fact that $\contatto\wedge(d\contatto)^{\wedge r}$ is a well-defined
section of $\Omega^{2r+1}_\sympy\tens[\sho]\sho_\sympy(r+1)$ for $0\leq r
\leq n$.)

There is an open covering $\shv=\{V_i\}_{i\in I}$ of 
$\sympy$ and  
contact embeddings $\chi_i \cl V_i \hookrightarrow P^*Y_i$. 
Up to refining the covering (we still denote it by $\shv$), the 
induced contact transformations 
$\chi_{ij}\cl \chi_i(V_{ij}) \to \chi_j(V_{ij})$ 
can be quantized to a $*$-preserving filtered $\C$-algebra isomorphism
$$
\Phi_{ij}\cl\chi_{ij}^{-1}(\Ev[Y_j]|_{\chi_j(V_{ij})})\to 
\Ev[Y_i]|_{\chi_i(V_{ij})}.
$$
The composition $\Phi_{ij}\Phi_{jk}\Phi_{ik}^{-1}$ is a 
$*$-preserving automorphism of $\Ev[Y_i]|_{\chi_i(V_{ijk})}$, 
and hence is equal to $\mathrm{Ad}(Q_{ijk})$ for a unique 
$Q_{ijk}\in\Ev[Y_i](V_{ijk})$ satisfying \eqref{eq:Qcond}. 
This proves the theorem below, thanks to Remark~\ref{re:aafaa}.

\begin{theorem} 
\label{th:Kq}
{\rm (cf~\cite[Theorem 2]{K_q})}
The triplet
\begin{equation}
\label{eq:AlgEE}
\EE[\sympy] = 
(\{\opb \chi_i \Ev[Y_i]|_{V_i}\}_{i\in I},
\{\opb\chi_i(\Phi_{ij})\}_{i,j\in I},\{\opb
\chi_i(Q_{ijk})\}_{i,j,k\in I})
\end{equation}
is a $\C$-algebroid descent datum over $\sympy$. 
In particular, there is an associated stack
$\stkMod(\EE[\sympy])$ on $\sympy$ locally equivalent to the stack 
$\stkMod(\Ev[Y_i])$ of microdifferential modules.
\end{theorem}

\subsection*{Good modules}

Any local notion, such as that of being coherent, simple, or regular, immediately extends to the category
$\catMod(\EE[\sympy])$. Here, we will discuss the non local notion of being good (cf~\cite{K_mu}).

Remark first that the construction in Theorem~\ref{th:Kq} also applies when 
replacing $\E[Y_i]$ with $\E[Y_i](0)$. One thus gets a $\C$-algebroid descent datum $\EE[\sympy](0)$ as well as a $\C$-linear functor $\EE[\sympy](0) \to \EE[\sympy]$. This induces a forgetful functor
$$
\catMod(\EE[\sympy]) \to[\for] \catMod(\EE[\sympy](0)).
$$
This functor locally admits a left adjoint, hence it has a left adjoint
$$
\catMod(\EE[\sympy](0)) \to[\exten] \catMod(\EE[\sympy]).
$$

\begin{definition}
A coherent module $\shm\in\catMod(\EE[\sympy])$ is good if for any relatively compact 
open subset $U\subset\sympy$ there exists $\shm_0\in\catMod(\EE[\sympy](0)|_U)$ such that 
$\shm|_U\simeq\exten(\shm_0)$.
\end{definition}

\subsection*{Quantization with parameters}

Assume now that the bundle $\sho_\sympy(1)$ is trivial, {\em i.e.}\ that
there exists a nowhere vanishing section
$\tau\in\sect(\sympy;\sho_\sympy(1))$. Consider an open covering
$\{V_i\}_{i\in I}$ of $\sympy$ and
contact embeddings $\chi_i \cl V_i \hookrightarrow P^*Y_i$ to which
is attached the $\C$-algebroid descent datum \eqref{eq:AlgEE}.

\begin{lemma}
Up to refining the covering, there exist
$*$-preserving QCT's
$\Phi_{ij}$ above $\chi_{ij}$, and sections
$T_i\in\sect(\chi_i(V_i);\Ev[Y_i](1))$, such that
$$
\sigma(T_i)\circ\chi_i = \tau,\quad T_i^* = -T_i, \quad
\Phi_{ij}(T_j) = T_i.
$$
\end{lemma}

\begin{proof}
Let $(t,x_1,\dots,x_n,u_1,\dots,u_n)$ be a local coordinate system
on $V_i$ such that the contact form is given by $dt+\sum_i u_idx_i$
and $\tau = \sigma(\partial_t)\circ\chi_i$. Then set $T_i =
\partial_t$.

Let $\Phi_{ij}^\dagger$ be a QCT above $\chi_{ij}$ such that
$\Phi_{ij}^\dagger(T_j) = T_i$. Then the proof goes as that of
\cite[Lemma~5.3~(iii)]{PS}.
Consider the QCT above the identity given by
$$
\Phi_{ij}^\ddagger \eqdot
\Phi_{ij}^{\dagger-1}\circ *\circ \Phi_{ij}^\dagger\circ *.
$$
There exists $Q\in\Ev[Y_j]$ such that
$\Phi_{ij}^\ddagger=\rmad(QQ^*)$.
Hence $\Phi_{ij}\eqdot \Phi_{ij}^\dagger \circ\rmad(Q)$
is $*$-preserving QCT above $\chi_{ij}$, and $\Phi_{ij}(T_j) = T_i$.
\end{proof}

Denote by $\Ev[Y_i,T_i]$ the subalgebra of $\Ev[Y_i]$ of operators
commuting with $T_i$.
As above, denote by $Q_{ijk}\in\Ev[Y_i](\chi_i(V_{ijk}))$
the unique operator
satisfying \eqref{eq:Qcond} such that
$\Phi_{ij}\Phi_{jk}\Phi_{ik}^{-1} = \rmad(Q_{ijk})$.
Since $\rmad(Q_{ijk})(T_i) = T_i$, it follows that $Q_{ijk}$
is a section of $\Ev[Y_i,T_i]$.

\begin{proposition}
\label{pr:EEtau}
Let $\sympy$ be a complex contact manifold.
Assume that there exists a nowhere vanishing section
$\tau\in\sect(\sympy;\sho_\sympy(1))$. Then
the triplet
\begin{equation}
\label{eq:AlgEEtau}
\EE[\sympy,\tau] =
(\{\opb \chi_i \Ev[Y_i,T_i]|_{\chi_i(V_i)}\}_{i\in I},
\{\opb\chi_i(\Phi_{ij})\}_{i,j\in I},\{\opb
\chi_i(Q_{ijk})\}_{i,j,k\in I})
\end{equation}
is a $\C$-algebroid descent datum on $\shv$.
In particular, there is an associated stack
$\stkMod(\EE[\sympy,\tau])$ on $\sympy$ locally equivalent
to the stack $\stkMod(\Ev[Y_i,T_i])$.
\end{proposition}

\section{Simple holonomic modules on contact manifolds}

Here, we give a proof of a result of Kashiwara~\cite{K_q} on the 
existence of twisted simple holonomic modules along smooth 
Lagrangian submanifolds (also called Legendrian in the
literature) of complex contact manifolds. 

Let $\sympy$ be a complex contact manifold. Recall from
Theorem~\ref{th:Kq} 
that there is an algebroid descent datum $\EE[\sympy]$ of
microdifferential 
operators on $\sympy$. 
Let $\Lambda\subset \sympy$ be a Lagrangian submanifold.
With notations as in Example~\ref{ex:Omega1/2}, consider the stack 
$\stkMod(\EE[\sympy]|_\Lambda \tens[\C] \C_{\sqrt{\Omega_\Lambda}})$ 
of twisted microdifferential modules on $\Lambda$. 

Let us recall Kashiwara's theorem announced in \cite{K_q}. 

\begin{theorem}\label{th:EEsimple}
Let $\sympy$ be a complex contact manifold and let  
$\Lambda\subset\sympy$ be a Lagrangian submanifold. 
There exists 
$\shl \in \catMod(\EE[\sympy]|_\Lambda \tens[\C]
\C_{\sqrt{\Omega_\Lambda}})$ 
which is simple along $\Lambda$.
\end{theorem}

\begin{remark}
It will follow from the proof that moreover $\shl$ is good.
\end{remark}

Here, we give a proof of this result using the notion of symbol 
for simple sections of holonomic modules recalled in Section~\ref{se:mumo}.

\begin{proof}
Let us denote for short by
$$
\EE[\sympy] = (\{\sha_i\}_{i\in I},\{f_{ij}\}_{i,j \in
I},\{a_{ijk}\}_{i,j,k \in I})
$$
the $\C$-algebroid descent datum \eqref{eq:AlgEE} attached to an open
covering $\sympy=\bigcup\nolimits_{i\in I}V_i$. 
Up to a refinement, we may assume that the $\C$-algebroid descent
datum $\C_{\sqrt{\Omega_\Lambda}}$ is attached to the same covering.
More precisely, using notations as in Example~\ref{ex:Omega1/2} for
$X=\Lambda$ and $U_i = \Lambda_i = \Lambda \cap V_i$, we assume that
there are nowhere vanishing sections $\omega_i\in\Omega_{\Lambda_i}$,
and $s_{ij} \in \OO^\times_{\Lambda_{ij}}$ 
with $\sqrt{\omega_j} = s_{ij} \sqrt{\omega_i}$, such that
$$
\C_{\sqrt{\Omega_X}} = ( \{ \C_{\Lambda_i} \}_{i \in I}, 
\{\id_{\C_{\Lambda_{ij}} } \}_{i,j \in I}, \{c_{ijk}\}_{i,j,k \in I}),
$$
where the $c_{ijk}$'s are defined by
$$
s_{ij}s_{jk} = c_{ijk} s_{ik}.
$$
Up to a further refinement of the covering, we may assume that 
there exist simple $\sha_i$-modules $\shl_i$ of order $0$ along 
$\Lambda_i$ and simple generators $u_i$ of $\shl_i$ such that 
$$
\sigma_{\Lambda_i}(u_i) = [\sqrt{\omega_i}] \in (
\sqrt{\Omega_\Lambda^\times}/\C_\Lambda^\times) |_{\Lambda_i}.
$$

Since $\shl_i$ and $\atw{f_{ji}}\shl_j$ are simple $\sha_i$-modules
along 
$\Lambda_{ij}$ of the same order, they are isomorphic, and one has
$$
\hom[\sha_i](\atw{f_{ji}}\shl_j,\shl_i) \simeq \C_{\Lambda_{ij}}.
$$ 
Let us describe explicitly such an isomorphism. 
For $\widetilde\xi_{ij}\in \hom[\sha_i](\atw{f_{ji}}\shl_j,\shl_i)$, let 
$\widetilde b_{ij}$ is a section of $\sha_i$ satisfying
$\widetilde\xi_{ij}(u_j) = \widetilde b_{ij} u_i$. 
One has 
\begin{align*}
\sigma_\Lambda(\widetilde\xi_{ij}(u_j)) 
&= \sigma_\Lambda(u_j) = [\sqrt{\omega_j}] = [s_{ij}\sqrt{\omega_i}],
\\
\sigma_\Lambda(\widetilde b_{ij}u_i) 
& = \sigma(\widetilde b_{ij})\sigma_\Lambda(u_i) = [\sigma(\widetilde
b_{ij})\sqrt{\omega_i}], 
\end{align*}
where the fourth equality follows from \eqref{eq:sigmaPu}. In particular,
$s_{ij}^{-1} \sigma(\widetilde b_{ij}) \in
\C^\times_{\Lambda_{ij}}$.
We can thus consider the isomorphism
\begin{align*}
\varphi_{ij}\cl\hom[\sha_i](\atw{f_{ji}}\shl_j,\shl_i) &\isoto
\C_{\Lambda_{ij}} \\
\widetilde\xi_{ij} &\mapsto  s_{ij}^{-1} \sigma(\widetilde b_{ij}).
\end{align*}
Set 
$$
\xi_{ij} = \varphi_{ij}^{-1}(1) \colon \atw{f_{ji}}\shl_j \isoto
\shl_i,
$$
and let $b_{ij}\in \sha_i$ satisfy $\xi_{ij}(u_j) = b_{ij} u_i$. As
$\varphi_{ij}(\xi_{ij}) = 1$, one has
\begin{equation}
\label{eq:sbij}
\sigma(b_{ij}) = s_{ij}.
\end{equation}
By Proposition~\ref{pr:twistedglue}, 
we get a twisted module
$$
\shl = (\shl_i,\xi_{ij})\in \catMod(\EE[\sympy]|_\Lambda \tens[\C]
\stks),
$$
where 
$$
\stks = ( \{ \C_{\Lambda_i} \}_{i \in I}, 
\{\id_{\C_{\Lambda_{ij}} } \}_{i,j \in I}, \{\widetilde
c_{ijk}\}_{i,j,k \in I}),
$$ 
the $\widetilde c_{ijk}$'s being given by \eqref{eq:cijk}. 
Since $\varphi_{ij}(\xi_{ij}) = 1 = \varphi_{jk}(\xi_{jk})$, by
\eqref{eq:phijikjki} we have
$$
\widetilde c_{ijk} =
\varphi_{ik}(\xi_{ij}(\atw{f_{ji}}\xi_{jk}(a_{kji}\cdot\scbul))).
$$
To show that $\stks = \C_{\sqrt{\Omega_\Lambda}}$, it thus remains to
prove that $\widetilde c_{ijk} = c_{ijk}$, i.e.~that
$$
s_{ij}s_{jk} = \widetilde c_{ijk} s_{ik}.
$$
One has 
\begin{equation*}
\begin{split}
\xi_{ij}(\atw{f_{ji}}\xi_{jk}(a_{kji} u_k) )
&=  f_{ij}f_{jk}(a_{kji}) \, \xi_{ij}  (\atw{f_{ji}}\xi_{jk} (u_k)) \\
&=  f_{ij}f_{jk}(a_{kji}) \, \xi_{ij} (b_{jk} u_j)   \\
&=  f_{ij}f_{jk}(a_{kji})  f_{ij}(b_{jk}) \, \xi_{ij}(u_j) \\
&=  f_{ij}f_{jk}(a_{kji})  f_{ij}(b_{jk})b_{ij} \, u_i .
\end{split}
\end{equation*}
Then
\begin{equation*}
\begin{split}
\widetilde c_{ijk} &= \varphi_{ik}( \xi_{ij}
(\atw{f_{ji}}\xi_{jk}(a_{kji} \cdot\scbul)) ) \\
&=  s_{ik}^{-1}  \sigma(f_{ij}f_{jk}(a_{kji})  f_{ij}(b_{jk})b_{ij})
\\
&=  s_{ik}^{-1}  \sigma(b_{jk})\sigma(b_{ij}) \\
&=  s_{ik}^{-1} s_{jk} s_{ij},
\end{split}
\end{equation*}
where the third equality follows from the fact that
$\sigma(a_{kji})=1$, and the fourth one from \eqref{eq:sbij}.
\end{proof}

\begin{definition}
\label{def:homsymmnf}
%(see~\cite{K3})
A complex homogeneous symplectic manifold $\stz =
(\stz,\simplettico,v)$ is a complex symplectic manifold
$(\stz,\simplettico)$ endowed with a vector field $v$ satisfying
$\Lie_v\simplettico=\simplettico$.
\end{definition}

Corollary~\ref{cor:kashi} below was announced in~\cite{K_q}. Although
we shall not use it, we give a proof for the reader's convenience.
Also note that our statement corrects that in loc.\ cit., following a
private communication with M.~Kashiwara.

Let $\sympy = (\sympy,\sho_\sympy(1),\contatto)$ be a contact manifold and
denote by 
$\gamma\cl\tw\sympy\to\sympy$ the total space of the
$\C^\times$-principal bundle associated with the dual of the line bundle
$\sho_{\sympy}(1)$. Then 
 $\tw\sympy$ is a  homogeneous symplectic manifold and there exists a
covering $\{V_i\}_{i\in I}$ of $\sympy$ by contact charts 
$\chi_i\colon V_i\hookrightarrow P^*Y_i$ 
such that, setting $\tw V_i = \opb\gamma V_i$, 
$\{\tw V_i\}_{i\in I}$ is a covering of $\tw \sympy$ and 
there are commutative diagrams 
\begin{equation*}
\xymatrix{
{\tw V_i\,\,}\ar[d]_{\gamma\vert_{\tw V_i}}\ar@{^(->}[r]^{\tw\chi_i}
                             &{\dot T^*Y_i}\ar[d]^{\gamma_i}\\
{V_i\,\,}\ar@{^(->}[r]^{\chi_i}&{P^*Y_i}.
}
\end{equation*}
where  $\gamma_i$ is the  
the projection $\dot T^*Y_i\to P^*Y_i$ and 
$\tw\chi_i\colon \tw V_i\hookrightarrow \dot T^*Y_i$ are 
homogeneous symplectic maps.

We denote by $\stkMod_{\text{loc-sys}}(\C_{\opb\gamma\Lambda})$ 
the full substack of
$\stkMod(\C_{\opb\gamma\Lambda})$ 
consisting of local systems, and by
$\stkMod_{\text{reg-}\Lambda}(\EE[\sympy]|_\Lambda)$ 
the full substack of $\stkMod(\EE[\sympy]|_\Lambda)$ consisting 
of modules
regular along $\Lambda$

\begin{corollary}\label{cor:kashi} \lp {\rm cf~\cite{K_q}}\rp\,
There is an equivalence of stacks
$$
\stkMod_{\text{reg-}\Lambda}(\EE[\sympy]|_\Lambda
\tens[\C]\C_{\sqrt{\Omega_\Lambda}}) \simeq
\oim\gamma \stkMod_{\text{loc-sys}}(\C_{\opb\gamma\Lambda}).
$$
\end{corollary}
\begin{proof}
The sheaf of rings $\E[Y]$ on $T^*Y$ is a subsheaf  of the
sheaf of rings $\ER[Y]$ of \cite{SKK}, 
this last sheaf being the microlocalization along the
diagonal of the sheaf $\OO_{Y\times Y}$ (up to a shift and
tensorizing by
holomorphic forms), see \cite[chapter 11]{KS1} for a detailed
construction.
One defines the sheaf $\ERv[Y]$ similarly as we have defined $\Ev[Y]$.

With notations as in \eqref{eq:AlgEE}, 
the isomorphisms $\Phi_{ij}$ are induced by sections of a
microdifferential 
module supported by the graph of $\chi_{ij}$, and hence extend to 
isomorphisms $\Phi_{ij}^\R$. Considering 
$Q_{ijk}$ as sections of $\ERv[Y_i]$, we  get a $\C$-algebroid 
descent datum
$$
\EER[\widetilde \sympy] = 
(\{\opb{\tw\chi_i}\ERv[Y_i]|_{\tw\chi_i(\tw V_i)}\}_{i\in I},
\{\opb{\tw\chi_i}(\Phi_{ij}^\R)\}_{i,j\in I},
\{\opb{\tw\chi_i}(Q_{ijk})\}_{i,j,k\in I}),
$$
and the associated stack $\stkMod(\EER[\tw \sympy])$ on 
$\tw\sympy$. The inclusion $\opb\gamma\Ev[Y] \subset \E[Y]^\R$ 
induces a functor of extension of scalars 
$\stkMod(\EE[\sympy]) \to \oim\gamma\stkMod(\EER[\widetilde \sympy])$. 
Let us denote by $\shm\mapsto\shm^\R$ this functor.

By Theorem~\ref{th:EEsimple} there is a simple system $\shl$ in 
$\catMod(\EE[\sympy]|_\Lambda \tens[\C] \C_{\sqrt{\Omega_\Lambda}})$. 
Consider the functor
$$
\Psi\cl\stkMod_{\text{reg-}\Lambda}(\EE[\sympy]|_\Lambda \tens[\C] 
\C_{\sqrt{\Omega_\Lambda}}) \to
\oim\gamma \stkMod_{\text{loc-sys}}(\C_{\opb\gamma\Lambda})
$$
given by $\shm \mapsto \hom[{\EER[\widetilde
\sympy]}](\shl^\R,\shm^\R)$.
The functor $\Psi$ is locally an equivalence by the results
of~\cite{KK}, and
hence is an equivalence.
\end{proof}

\begin{example}
Let $\sympy = P^*\C$, $\Lambda = P^*_{\{0\}}\C = \{\rmpt\}$,
$\gamma^{-1}\Lambda = \C^\times$. In this case, simple objects of
rank one are classified by $\C^\times\simeq\C/\Z$.
\end{example}

\section{WKB-modules}

The relationship between microdifferential operators on a 
complex contact manifold and
WKB-differential operators on a complex symplectic manifold is
classic, and is discussed {\em e.g.}, in \cite{AKKT, Ph} in the case
of
 cotangent bundles. This study (including the analysis of the 
action of quantized contact transformations) is systematically 
performed in~\cite{PS}, and here we follow their presentation. 

Let $X$ be a complex manifold, $t\in\C$ the coordinate, and set
$$
\E[X\times\C, \hat t] = \E[X\times\C, \partial_t] = \{P\in \E[X\times\C] ;\
[P,\partial/\partial_t] = 0\}.
$$
Set $\dot P^*(X\times\C) = \{\tau\neq 0\}$, and consider the
projection 
\begin{equation}
\label{eq:projrho}
\rho\cl \dot P^*(X\times\C) \to T^*X,
\end{equation}
given in local coordinates by $\rho(x,t;\xi,\tau) = (x;\xi/\tau)$.
The ring of WKB-operators on $T^* X$ is defined by
$$
\W[X] := \oim\rho(\E[X\times\C, \hat t]).
$$
We similarly set $\Wv[X] := \oim\rho(\Ev[X\times\C, \hat t])$.
In a local symplectic coordinate system $(x;\xi)$ on $T^*X$, a
section $P\in \W[X](U)$ is written as a formal series
$$
P = \sum_{j\leq m}p_j(x;\xi)\tau^j, \quad p_j\in\OO_{T^*X}(U),\quad
m\in\Z,
$$
with the condition that for any compact subset $K$ of $V$ there
exists a constant
$C_K>0$ 
such that $\sup\limits_{K}\vert p_{j}\vert \leq C_K^{-j}(-j)!$ for
all $j<0$.

One sets 
\begin{equation}
\field \eqdot \W[\rmptt].
\end{equation}
Hence, an element  $a\in\field$ is written as a formal series
$$
a = \sum_{j\leq m}a_j\tau^j, \quad a_j\in\C,\quad  m\in\Z,
$$
with the condition that there exist 
$C>0$ with $\vert a_{j}\vert \leq C^{-j}(-j)!$ for all $j<0$.

Note that $\W[X]$ is $\Z$-filtered $\field$-central algebra,
and the principal symbol map
$
\sigma_m \cl \W[X](m) \to \OO_{T^*X}
$
induces an isomorphism of graded algebras
$\gr\W[X] \to \OO_{T^*X}[\tau^{-1},\tau]$. Note also that
$\opb{\pi}\shd_X$ is a subring of $\W[X]$.

We can now mimic Definition \ref{def:simpleEmod} for $\W[X]$-modules.

\begin{definition}\label{def:simpleWmod}
Let $\Lambda$ be a smooth Lagrangian submanifold of $T^*X$. 
Let $\shm$ be a coherent $\W[X]$-module supported by $\Lambda$.  
One says that $\shm$ is regular (resp.~simple) along
$\Lambda$ if there locally  exists a coherent sub-$\W[X] (0)$-module
$\shm_0$
of $\shm$ which generates it over $\W[X] $, and such that
$\shm_0/\W[X] (-1)\shm_0$ is an $\OO_\Lambda$-module (resp.~a locally
free $\OO_\Lambda$-module of rank one).  
\end{definition}

Note that, as follows e.g.~from Corollary~\ref{co:Wreg} below, if
$\shm$ is regular then it is locally a finite direct sum of simple
modules.

\begin{notation}
Let $X$ be  a complex manifold. 
We denote by $\OO^\tau_X$ the simple
$\W[X]$-module along the zero-section $T^*_XX$ defined 
by $\OO^\tau_X = \W[X]/\shi$, where $\shi$ is the left ideal 
generated by the vector fields on $X$. 
\end{notation}

\begin{proposition}\label{pr:locSimp}
\begin{enumerate}[{\rm(i)}]
\item
Any two $\W[X]$-modules simple along $\Lambda$ are locally
isomorphic. In particular, any simple module along $T^*_XX$ is 
locally isomorphic to $\OO^\tau_X$.
\item
If $\shm$, $\shn$ are simple $\W[X]$-modules along
$\Lambda$, 
then
$\rhom[{\W[X]}](\shm, \shn)$ is a $\field$-local system of rank one
on 
$\Lambda$.
\end{enumerate}
\end{proposition}

\begin{proof}
Since both statements are local on $\Lambda$, 
we may assume that $X$ is endowed with a local coordinate system 
$(x_1,\dots,x_n)$ and that $\Lambda$ is the zero-section $T^*_XX$ of
$T^*X$. 

\noindent
(i) Any $\W[X]$-module simple along $T^*_XX$ is
locally isomorphic to $\sho_X^\tau$. Indeed, the proof 
 of the theorem (due to \cite{SKK})
 which asserts that if 
 $Y$ is a complex manifold, then 
simple $\E[Y]$-modules along smooth regular involutive
submanifolds of $P^*Y$ are locally isomorphic applies when replacing
$\E[X\times\C]$ with $\E[X\times\C, \hat t]$ (for an exposition, 
see \cite[Ch 1, Th 6.2.1]{Sc}). 

\noindent
(ii) By (i) we may assume that $\shm=\shn=\sho_X^\tau$. 
The result easily follows, representing 
$\rhom[{\W[X]}](\sho_X^\tau, \sho_X^\tau)$ by the Koszul complex 
$K^\bullet(\sho_X^\tau,(\partial_{x_1},\dots,\partial_{x_n}))$
associated with the sequence $(\partial_{x_1},\dots,\partial_{x_n})$
acting on $\sho_X^\tau$.
\end{proof}

Recall the projection
$\rho$ in \eqref{eq:projrho}
and note that $\opb\rho\Lambda$ is an involutive submanifold 
of $\dot P^*(X\times\C)$.

\begin{proposition}
\label{pr:Lambdac}
Let $\Lambda$ be a smooth Lagrangian submanifold of $T^*X$.
\begin{enumerate}[{\rm(i)}]
\item
Locally, there exists a Lagrangian submanifold
$\Lambda^0\subset\rho^{-1}(\Lambda)$ on which $\rho$
induces an isomorphism $\Lambda^0\simeq\Lambda$.
\item
Let $\shl$ be a simple $\W[X]$-module along $\Lambda$. Then, locally
there exists 
a simple $\E[X\times\C]$-module $\shl^0$ along $\Lambda^0$ such that
$\shl \simeq \for(\shl^0)$, where 
$\for$ denotes the forgetful functor $\catMod(\E[X\times\C])
\to\catMod(\E[X\times\C,\hat t])$.
\end{enumerate}
\end{proposition}

\begin{proof}
Since the problem is local, we may assume that $\Lambda=T^*_XX$
in a system of local symplectic coordinates $(x;u)\in T^*X$. In 
the corresponding system of homogeneous coordinates 
$(x,t;\xi,\tau)\in P^*(X\times\C)$, one has
$\opb\rho\Lambda=\{\xi=0\}$. 
This set is foliated by the Lagrangian 
submanifolds $\Lambda^c=\{\xi=0,\,t=c\}$, 
for $c\in\C$. For $Z\subset X$ a closed submanifold, 
denote by $\shc_{Z|X}$ the sheaf of finite order holomorphic
microfunctions on
$P^*_ZX$. 
A simple $\E[X]$-module along $\Lambda^0$ is given by 
$\shl^0 = \shc_{X\times\{0\}|X\times\C}$. 
One then immediately checks that $\shl^0\simeq\OO^\tau_X$.
\end{proof}

\section{Contactification of symplectic manifolds}

In this section, we recall some well-known facts from the specialists
on contact and symplectic geometry.

\begin{definition}\label{de:contactf}
A contactification of a symplectic manifold
$(\sympx,\simplettico)$ is a complex contact manifold
$(\sympy,\sho_\sympy(1),\contatto)$,
and a morphism of complex manifolds
$\rho\colon\sympy\to\sympx$ such that
the following conditions are satisfied:
\begin{enumerate}[{\rm(a)}]
\item
the line bundle $\sho_\sympy(1)$
has a nowhere vanishing section $\tau$,
\item there exists an open covering $\sympx = \bigcup\nolimits_{i\in I}
U_i$,
holomorphic functions $t_i$ on $\opb{\rho}(U_i)$  and
primitives $\sigma_i\in\Omega^1_{U_i}$ of
$\simplettico|_{U_i}$ such that
$\chi_i\eqdot (\rho , t_i)$ gives an isomorphism
$\chi_i\cl\opb{\rho}(U_i)\isoto U_i\times\C$, and
$dt_i+\rho^*(\sigma_i) = (\contatto/\tau)|_{\opb\rho(U_i)}$.
\end{enumerate}
\end{definition}

\begin{lemma}
\label{le:affine}
A contactification $\rho\colon\sympy\to\sympx$ has a structure of
$\C$-principal bundle.
\end{lemma}

\begin{proof}
With notations as in Definition~\ref{de:contactf}, set $\chi_{ij} =
\chi_j\chi_i^{-1}|_{U_{ij}} \colon U_{ij}\times\C \to U_{ij}\times\C$. We have
$\chi_{ij} = (\id_{U_{ij}}, \theta_{ij})$ for some transition functions
$\theta_{ij}\colon U_{ij}\times\C\to\C$. Since $d(t_i - t_j) =
\rho^*(\sigma_i-\sigma_j)$, up
to shrinking the covering we may assume that $\sigma_i - \sigma_j =
df_{ij}$ for some functions $f_{ij}\in\sho_{U_{ij}}$. The transition functions 
are the translations $\theta_{ij}(p,t) = t +f_{ij}(p)+c_{ij}$ for some $c_{ij}\in\C$.
\end{proof}

\begin{lemma}
\label{le:zeroclass}
A symplectic manifold $(\sympx,\simplettico)$ admits a contactification
if and only if the de~Rham cohomology class $[\simplettico]\in
H^2(\sympx;\C_\sympx)$ vanishes.
\end{lemma}

\begin{proof}
(a)
Let $\sympx = \bigcup\nolimits_{i\in I} U_i$ be a covering such that
$\simplettico|_{U_i}$ has a primitive $\sigma_i\in\Omega^1_{U_i}$. Up
to shrinking the covering we may assume that $\sigma_i - \sigma_j =
df_{ij}$ for some functions $f_{ij}\in\sho_{U_{ij}}$. Then $d(f_{ij}
+ f_{jk} - f_{ik}) = 0$, and $\C_{U_{ijk}} \owns c_{ijk} = f_{ij} +
f_{jk} - f_{ik}$ is a Cech cocycle representing $[\simplettico]\in
H^2(\sympx;\C_\sympx)$. Since $[\simplettico] = 0$, we have $c_{ijk} =
c_{ij} + c_{jk} - c_{ik}$. Setting $\widetilde f_{ij} = f_{ij} -
c_{ij}$, we have $\widetilde f_{ij} + \widetilde f_{jk} - \widetilde
f_{ik} = 0$. 
Endow $U_i\times \C$ with the contact form $dt+\sigma_i$.
A contactification of $\sympx$ is thus given by the $\C$-principal
bundle $\sympy\to\sympx$ with local charts $U_i\times \C$, and transition
functions over $U_{ij}$ given by $(p,t)\mapsto (p,t+\widetilde
f_{ij}(p))$.

\smallskip\noindent(b)
Let $\rho\colon\sympy\to\sympx$ be a contactification, and use notations
as in Definition~\ref{de:contactf}.
We have $d(\contatto/\tau) = \rho^*\simplettico$ and hence
$[\rho^*\simplettico]=0$ in $H^2(\sympy;\C_\sympy)$.
Then $[\simplettico]=0$ in $H^2(\sympx;\C_\sympx) \isoto
H^2(\sympy;\C_\sympy)$.
\end{proof}

\begin{lemma}
\label{le:contlag}
Let $(\sympx,\simplettico)$ be a symplectic manifold and $\Lambda$ a
Lagrangian submanifold. After replacing $\sympx$ with a neighborhood of
$\Lambda$ there exists a contactification $\rho\colon \sympy \to \sympx$
and a Lagrangian submanifold $\Lambda^0 \subset \sympy$ on which $\rho$
induces an isomorphism $\Lambda^0 \simeq \Lambda$.
\end{lemma}

\begin{proof}
The restriction map $H^2(\sympx,\C_\sympx) \to H^2(\Lambda;\C_\Lambda)$
is given by $[\simplettico] \mapsto [j^*\simplettico]$, where
$j\colon \Lambda \hookrightarrow \sympx$ is the embedding. Since
$\simplettico|_\Lambda =0$, there exists an open neighborhood
$U\supset \Lambda$ such that $[\simplettico|_U] = 0$. Thus, up to
replacing $\sympx$ with $U$, we can assume that $[\simplettico] = 0$.

Let us adapt the arguments in part (a) of the proof of
Lemma~\ref{le:zeroclass} above.
Let $\sympx = \bigcup\nolimits_{i\in I} U_i$ be a covering such that
$\simplettico|_{U_i}$ has a primitive $\sigma_i\in\Omega^1_{U_i}$. We
may assume that $\sigma_i|_\Lambda = 0$.
Let $f_{ij}\in\sho_{U_{ij}}$ satisfy $\sigma_i - \sigma_j = df_{ij}$.
We may assume moreover that $f_{ij}|_\Lambda = 0$. Since $f_{ij} +
f_{jk} - f_{ik}$ is locally constant on $\Lambda$, it must vanish on $\Lambda$. 
Thus, the $\C$-principal bundle $\rho\colon \sympy\to\sympx$ is described by the
local charts $U_i\times \C$, with transition functions 
$(p,t)\mapsto(p,t+f_{ij}(p))$ and contact form $dt+\sigma_i$.
We define $\Lambda^0$ by $\Lambda^0|_{U_i\times\C} =
\Lambda|_{U_i}\times\{0\}$.
\end{proof}

\section{Simple holonomic modules on symplectic
manifolds}\label{se:sympl}

Let us start by recalling the construction of \cite{PS} of
the stack of WKB-modules on a symplectic manifold,
in the special case where there exists a contactification.

Let $\rho\colon\sympy\to\sympx$ be a contactification of a complex
symplectic manifold. Denote by
$$
\EE[\sympy,\tau] =
(\{\sha_i\}_{i\in I},\{f_{ij}\}_{i,j \in I},\{a_{ijk}\}_{i,j,k \in I})
$$
the $\C$-algebroid descent datum on $\sympy$
given by Proposition~\ref{pr:EEtau}.
Consider the $\C$-algebroid descent datum on $\sympx$
$$
\WW[\sympx] =
(\{\oim\rho\sha_i\}_{i\in I},\{\oim\rho f_{ij}\}_{i,j \in I},
\{\oim\rho a_{ijk}\}_{i,j,k \in I}).
$$
The  stack $\stkMod(\WW[\sympx])$  associated with this $\C$-algebroid
descent datum
is the stack of WKB-modules in~\cite{PS}.

One defines the notion of good $\WW[\sympx]$-module similarly as for $\EE[\sympy]$-modules.

Here, we prove the existence of twisted simple holonomic modules
along 
smooth Lagrangian submanifolds of complex symplectic manifolds.

\begin{theorem}\label{th:WWsimple}
Let $\sympx$ be a complex symplectic manifold and let $\Lambda\subset
\sympx$ be a Lagrangian submanifold.
There exists a module
$\shl \in \catMod(\WW[\sympx]|_\Lambda \tens[\C]
\C_{\sqrt{\Omega_\Lambda}})$ 
which is simple along $\Lambda$.
Moreover, $\shl$ is good.
\end{theorem}

\begin{proof}
By Lemma~\ref{le:zeroclass} there exist a contactification 
$\rho\colon\sympy\to\sympx$ and a Lagrangian submanifold 
$\Lambda^0\subset \sympy$ on which $\rho$ induces an 
isomorphism $\Lambda^0\simeq\Lambda$. By Theorem~\ref{th:EEsimple} 
there exists $\shl^0 \in \catMod(\EE[\sympy]|_{\Lambda^0} \tens[\C]
\C_{\sqrt{\Omega_{\Lambda^0}}})$ 
which is simple along $\Lambda^0$. As in Proposition~\ref{pr:Lambdac}
we then set
$\shl = \for(\shl^0)$,
where
\begin{align*}
\for \colon \catMod(\EE[\sympy]|_{\Lambda^0} \tens[\C]
\C_{\sqrt{\Omega_{\Lambda^0}}}) &\to \catMod(\EE[\sympy,\hat
t]|_{\Lambda^0} \tens[\C]
\C_{\sqrt{\Omega_{\Lambda^0}}}) \\
&\simeq \catMod(\WW[\sympx]|_\Lambda \tens[\C]
\C_{\sqrt{\Omega_\Lambda}})
\end{align*}
is the natural forgetful functor.
\end{proof}

\begin{corollary}
\label{co:Wreg}
There is a $\field$-equivalence of stacks
\begin{equation}
\label{eq:WRS}
\stkMod_{\text{reg-}\Lambda}(\WW[\sympx]|_\Lambda \tens[\C] 
\C_{\sqrt{\Omega_\Lambda}}) \simeq
\stkMod_{\text{loc-sys}}(\field_\Lambda).
\end{equation}
\end{corollary}

\begin{proof}
By Theorem~\ref{th:WWsimple}, there exists a simple module 
$\shl$ in
$\catMod(\WW[\sympx]|_\Lambda\tens[\C]\C_{\sqrt{\Omega_\Lambda}})$. 
By Proposition~\ref{pr:locSimp}, a functor \eqref{eq:WRS} 
is given by $\shm \mapsto \hom[{\WW[\sympx]}](\shl,\shm)$.  Proving
that it
is an equivalence is a local problem, and so we may assume
$\sympx=T^*X$,
$\Lambda = T^*_XX$.
Then a quasi inverse is given by $F \mapsto F\tens[\field]\OO^\tau_X$.
\end{proof}

\begin{remark}
Note that in the real setting, a link between
simple holonomic modules on real  Lagrangian submanifolds and Fourier
distributions, Maslov index, etc.\  is investigated in \cite{NT}.
\end{remark}

\providecommand{\bysame}{\leavevmode\hbox to3em{\hrulefill}\thinspace}

\end{document}